\magnification=\magstep1
\hsize=16.5 true cm 
\vsize=23.6 true cm
\font\bff=cmbx10 scaled \magstep1
\font\bfff=cmbx10 scaled \magstep2
\font\bffg=cmbx10 scaled \magstep3

\font\smc=cmcsc10 
\parindent0cm
\overfullrule=0cm
\def\cl{\centerline}           %
\def\rl{\rightline}            %
\def\bp{\bigskip}              %
\def\mp{\medskip}              %
\def\sp{\smallskip}            %
\font\boldmas=msbm10           %
\def\Bbb#1{\hbox{\boldmas #1}} %
\def\Z{\Bbb Z}                 %
\def\Q{\Bbb Q}                 %
\def\R{\Bbb R}                 %
                 %
\def\({\left(}              %
\def\){\right)}             %
\def\[{\left[}              %
\def\]{\right]}             %

\def\={\;=\;}               %
\def\+{\,+\,}               %
\def\<{\;<\;}               %
\def\>{\;>\;}               %

\def\D{{\cal D}}
\def\G{{\cal G}}

\def\Y{{\cal R}}

\def\deg{\hbox{\rm deg}\,}
\cl{\bffg On the distribution of reducible polynomials}
\bp
\cl{\bfff Gerald Kuba}
\bp\mp
{\bff 1. Introduction and statement of results}
\mp
For a fixed integer $\,n\geq 2\,$ and a real parameter
$\,t\geq 1\,$ we consider all polynomials
\hbox{$p(X)\,=\,a_nX^n+a_{n-1}X^{n-1}+\cdots+a_2X^2+a_1X+a_0\;$}
with coefficients $\,a_i\in\Z\,$ such that $\,a_n\not=0\,$
and $\;H(p)\leq t\;$ 
where $\;H(p)\,:=\,\max\,\{\,|a_i|\;|\;i=0,1,...,n\,\,\}\;$
is the {\it height} of $\,p(X)\,$.
Of course, the total number of all these polynomials 
equals $\;[2t]\cdot[2t+1]^n\,\asymp\,t^{n+1}\;$
where $\,[\;]\,$ are the Gauss brackets
and $\;A\asymp B\;$ means $\;B\ll A\ll B\,$.
Let $\,\Y_n(t)\,$ denote the set of all 
polynomials $\,p(X)\,$ over $\Z$ 
with degree $\;n\geq 2\;$ and height $\;\leq t\;$
which are {\it reducibel} over $\Q$.
Note that $\,p(X)\,$ is reducible in the ring $\,\Q[X]\,$ if and only 
if $\,p(X)\,$ can be written as a product of two
polynomials in the ring  $\,\Z[X]\,$ both of less degree than $\,p(X)\,$.
\sp
In the famous exercise book of {\smc P¢lya} and {\smc Szeg{\"o}}
[3, Ex.266] one can find the estimate
\sp
\cl{$|\Y_n(t)|\,=\,O\big(t^n(\log t)^2\big)\qquad(t\to\infty)\,$.}
\mp
By the method used in {\smc D\"orge} [2]
this estimate can be improved 
to $\;|\Y_n(t)|=O\big(t^n(\log t)\big)\;$ which still
is not best possible when $\;n\geq 3\,$.
Indeed, the {\it true order of magnitude}
of the lattice points counting function 
$\;t\mapsto|\Y_n(t)|\;$ reads as follows.
\mp
{\bf Theorem 1.} {\it  For every integer $\,n\geq 3\,$
there is a constant $\;C_n>0\;$ such that 
\sp
\cl{$t^n\;\leq\;|\Y_n(t)|\;\leq\;C_n\cdot t^n$}
\sp
for all $\,t\geq 1\,$.}
\mp\sp
{\bf Theorem 2.} {\it As $\;t\to\infty\,$, 
$\;\;t^2\log t\;\ll\;|\Y_2(t)|\;\ll\;t^2\log t\;.$}
\mp\sp
There is a natural generalization of Theorem 2. 
Let $\,R_n^s(t)\,$ denote the total number 
of all polynomials $\,p(X)\,$ over $\Z$ 
with degree $\,n\,$ and height $\;\leq t\;$
such that $\,p(X)\,$ splits completely into linear factors in 
the ring $\,\Q[X]\,$ or, equivalently, in the ring $\,\Z[X]\,$.
Naturally, $\;R_2^s(t)=|\Y_2(t)|\,$, so that by Theorem 2 we have 
$\;\;R_2^s(t)\,\asymp\,t^2\log t\;\;$
as $\;t\to\infty\,$.  Now in general the following
estimation holds.
\mp
{\bf Theorem 3.} {\it For every fixed $\;n\geq 2\,$, 
$\;\;t^2(\log t)^{n-1}\,\ll\,R_n^s(t)\,\ll\,t^2(\log t)^{n-1}\;\;
(t\to\infty)\;.$}
\mp
Certainly, Theorem 3 is also true in the trivial case
$\;n=1\;$ where $\;R_1^s(t)\,=\,[2t]\cdot[2t+1]\;$
for every $\;t\geq 1\,$. (Of course, in general the case
$\,n=1\,$ is of no interest since $\;\Y_1(t)=\emptyset\,$.)
\mp
Theorem 3 demonstrates that the totally splitting polynomials
contribute only very little to the total number of
all reducible polynomials of fixed degree $\;\geq 3\,$ and
bounded height. On the other hand there is a special subclass
of $\,\Y_n(t)\,$ whose contribution 
to $\,|\Y_n(t)|\,$ is absolutely dominating.
This class lies on top of a hierarchy
of pairwise disjoint subclasses of $\,\Y_n(t)\,$.
For $\;{n\over 2}<k<n\;$ fixed and arbitrary $\;t\geq 1\;$
let $\;\Y_{k,n}(t)\subset\Y_{n}(t)\;$
such that $\;p(X)\in\Y_{k,n}(t)\;$ if and only if
$\,p(X)\,$ has an irreducible factor in $\,\Z[X]\,$
of degree $\,k\,$. The following theorem
shows that $\;\Y_{n-1,n}(t)\;$ is the mentioned top class.
\mp
{\bf Theorem 4.} {\it For $\;1<{n\over 2}<k<n\;$ fixed we
have $\;\;t^{k+1}\,\ll\,|\Y_{k,n}(t)|\,\ll\,t^{k+1}\;\;\;(t\to\infty)\,$.
\sp           
Specifically, 
$\;\;|\Y_{n-1,n}(t)|\,\asymp\,|\Y_{n}(t)|\,\asymp\,t^n\;\;\;(t\to\infty)
\;\;$ for every $\;n\geq 3\,$.
\sp
Moreover,
$\;\;|\,\Y_{n}(t)\setminus\Y_{n-1,n}(t)\,|\,\ll\,
t^{n-1}(\log t)^2\;\;\;(t\to\infty)\;\;$
for every $\;n\geq 3\;$
and the factor $\,(\log t)^2\,$ can be omitted if and only if $\;n\geq 4\,$.} 
\eject
\bp
{\bff 2. Preparation of the proofs}
\mp
Since $\;e^x>2^{x-2}\sqrt{x+1}\;$ for $\,x\geq 0\,$,
as an immediate consequence of [4] Theorem 4.2.2 we obtain
\mp\sp
{\bf Lemma 1.} {\it If $\;p,q\;$ are polynomials over $\,\Z\,$
with positive degrees $\;\deg p\;$ and $\;\deg q\;$ 
such that $\;n\,=\,\deg(pq)\,=\,\deg p\,+\,\deg q\;$
then $\;\;e^{-n}\,H(p)\,H(q)\;\leq\;H(pq)
\;\leq\;n\,H(p)\,H(q)\;.$} 
\bp
For $\;T\geq 1\;$ consider the hyperbola triangle
\mp
\cl{$\;\D(T)\;:=\;\{\,(x,y)\in\R^2\;|\;x,y\geq 1\;\;\land\;\;xy\leq T\,\}\;$}
\mp
and define the integral

$$I(T;a,b)\;\;:=\;\;
\int\!\!\!\!\int\limits_{\!\!\!\!\D(T)} x^ay^b\,d(x,y)$$

with real exponents $\;a,b\geq 0\,$. We compute

$$I(T;a,b)\;\;=\;\;{1\over (a+1)(b+1)}\,+\,{1\over a-b}
\Big({T^{a+1}\over a+1}-{T^{b+1}\over b+1}\Big)\qquad\hbox{\it when}\quad 
a\not=b$$

and 

$$I(T;c,c)\;\;=\;\;
{T^{c+1}\log T\over c+1}\,-\,{T^{c+1}-1\over (c+1)^2}
\qquad\hbox{\it for}\quad c\geq 0$$
\sp
and hence we obtain
\mp\sp
{\bf Lemma 2.} {\it For $\;a,b\geq 0\;$ fixed we have
\sp
\cl{$\;\;I(T;a,b)\;\asymp\;T^{1+\max\{a,b\}}\,(\log T)^\nu\;\;
\;(\,T\to\infty\,)\;\;$}
\sp
with $\;\nu=0\;$ when  $\;a\not= b\;$ and $\;\nu=1\;$ when  $\;a= b\;$.}
\bp
As usual, let $\;\varphi(\cdot)\;$ denote the Euler totient function.
We will use the following
well-known result due to {\smc Mertens} (Theorem 22 in [1]).
\mp\sp
{\bf Lemma 3.} {\it As $\;t\to\infty\,$, 
$\;\;\sum\limits_{m\leq t}\varphi(m)\;=\;{3\over \pi^2}t^2\,+\,O(t\log t)\;$.}
\bp
Further we will need the following lemma which 
immediately follows from Lemma 3 via partial summation.
\mp\sp
{\bf Lemma 4.} {\it As $\;t\to\infty\,$, 
\sp
\cl{$\sum\limits_{m\leq t}\varphi(m)\cdot m^{-2}\;\,\asymp\,\;\log t$}
\sp
and
\sp
\cl{$\sum\limits_{m\leq t}\varphi(m)\cdot m^\alpha\;\,\asymp\;\,
t^{\max\{0,\alpha+2\}}$}
\sp 
for every real $\;\alpha\,\not=\,-2\,$.
(The two $\asymp$-constants depend only on $\,\alpha\,$.)}
\vfill\eject
\bp
{\bff 3. Proof of Theorem 1}
\mp
The lower bound in Theorem 1 is trivial
since there are $\;[2t]\cdot[2t+1]^{n-1}\;$
polynomials $\,p(X)\,$ over $\,\Z\,$
with degree $\,n\,$ and height $\;\leq t\;$ such that $\;p(0)=0\,$. 
\mp
Let $\;{\cal P}_n(t)\;$ denote the set of all pairs $\,(p,q)\,$
of non-constant polynomials over $\Z$ such that $\;\deg p\,+\,\deg q\,=\,n\;$
and $\;H(pq)\leq t\,$. Then we obviously have
$\;|{\cal P}_n(t)|\geq |{\cal R}_n(t)|\;$ for all $\;t\geq 1\,$.
In view of Lemma 1 the set 
\sp
\cl{$\;{\cal P}_n^*(t)\;:=\big\{\,(p,q)\in(\Z[X]\setminus\Z)^2\;\;\big|\;\;
\deg p+\deg q\,=\,n\;\;\land\;\;H(p)\cdot H(q)\;\leq\;e^{n}t\,\big\}$}
\sp 
contains the set $\;{\cal P}_n(t)\;$ and thus we have
the estimate $\;|{\cal P}_n^*(t)|\geq |{\cal R}_n(t)|\;$ for all $\;t\geq 1\,$.
In order to prove Theorem 1 we show
\sp
(3.1)\qquad\qquad$|{\cal P}_n^*(t)|\;\ll\;t^n\;\;(t\to\infty)$
\sp
for every $\;n\geq 3\,$.
\mp
For abbreviation we set $\;T\,=\,e^{n}t\,$.
Since obviously
\sp
\cl{$|\{\,p\in\Z[X]\;\;|\;\;\deg p\,=\,k\;\;\land\;\; H(p)=h\,\}|
\;\leq\;2\cdot (2h+1)^k\cdot (k+1)$}
\sp
we obtain (3.1) by showing 
$$\sum_{k=1}^{n-1}\sum_{(x,y)\in \G(T)}
2(2x+1)^k(k+1)\cdot 2(2y+1)^{n-k}(n-k+1)
\;\;\ll\;\;T^n\;\;(T\to\infty)$$
where $\;\G(T)\;:=\;\D(T)\cap\Z^2\;$ with
$\;\D(T)\;=\;\{\,(x,y)\in\R^2\;|\;x,y\geq 1\;\;\land\;\;xy\leq T\,\}\,$.
Thus it is enough to verify
$$\sum_{(x,y)\in \G(T)}x^ky^{n-k}
\;\;\ll\;\;T^n\;\;\;(T\to\infty)\leqno(3.2)$$
for $\;1\leq k<n\;$ and $\;n\geq 3\;$ fixed.
\mp
Now, (3.2) is true because by Lemma 2 for the corresponding integral we have
\sp
\cl{$\;I(T,k,n-k)\,\ll\,T^n\;$}
\sp
provided that $\;n\geq 3\,$.
(Clearly, the difference between the sum in (3.2)
and $\;I(T,k,n-k)\;$ is $\;\ll\,T^n\;$ as $\;T\to\infty\,$.)
Additionally, $\;I(T;1,1)\;\ll\;T^{2}\log T\;$ 
yields $\;\;|{\cal P}_n^*(t)|\;\ll\;t^n\log t\;\;$ 
in the exceptional case $\,n=2\,$
and hence we also obtain the upper bound in Theorem 2.
\bp\mp
{\bff 4. Proof of Theorem 2}
\mp
It remains to verify the lower bound in Theorem 2.
As usual, we call a linear polynomial $\;aX+b\;$
over $\Z$ {\it primitive} when $\,a,b\,$ are coprime. 
Then for every quadratic polynomial $\,q(X)\,$ over $\Z$
which splits over $\Q$ there exists
one and only one set $\,\{f,g\}\,$ 
of primitive linear
polynomials $\,f(X),g(X)\,$ such that $\;f(X)g(X)\;$ divides $\;q(X)\;$      
in the ring $\,\Z[X]\,$. 
Let $\;Q(t)\;$ denote the number of all
quadratic polynomials over $\Z$ with height $\;\leq t\;$
which split over $\Z$ into two primitive linear factors.
Then we have $\;Q(t)\leq |{\cal R}_2(t)|\;$ and 
$\;2\cdot Q(t)\;$ is not smaller than the
cardinality of the set
\sp
\cl{$\big\{\,(f,g)\in\Z[X]^2\;\;\big|\;\;
\deg f\,=\,\deg g\,=\,1\;\;\land\;\;
(\,f,g\;\hbox{\it primitive}\,)\;\;\land\;\;
H(fg)\;\leq\;t\,\big\}\;$,}
\sp 
which contains the set
\sp
\cl{$\big\{\,(f,g)\in\Z[X]^2\;\;\big|\;\;
\deg f\,=\,\deg g\,=\,1\;\;\land\;\;
(\,f,g\;\hbox{\it primitive}\,)\;\;\land\;\;
H(f)\cdot H(g)\;\leq\;{1\over 2}t\,\big\}$}
\mp 
in view of Lemma 1. 
\sp
Further, the total number of all primitive
linear polynomials in $\,\Z[X]\,$
with constant height $\;h\geq 2\;$ 
clearly equals $\;8\cdot \varphi(h)\,$.
The number is equal to 6 when $\;h=1\,$.
Therefore with the new parameter
$\;T={1\over 2}t\;$ we have
$$|{\cal R}_2(t)|\;\;\geq\;\;
18\,\cdot\!\!\! \sum_{(x,y)\in \G(T)}\varphi(x)\,\varphi(y)$$
and the proof of Theorem 2 is finished by showing
$$T^2\log T\;\;\ll\sum_{(x,y)\in \G(T)}\varphi(x)\,\varphi(y)
\;\;\;(T\to\infty)\;.\leqno(4.1)$$

Note that (4.1) would 
immediately follow from $\;m\,\ll\,\varphi(m)\;$
and the fact that the sum in (3.2) with $\;k=1\;$ and $\;n=2\;$
is $\;\gg\,T^2\log T\,$, but of course $\;m\,\ll\,\varphi(m)\;$ is false
although $\;m^{1-\epsilon}\,\ll\,\varphi(m)\;$
is true for every $\,\epsilon>0\,$.
\mp                                           
Nevertheless we will reach our goal by applying Lemma 3.
As a consequence of Lemma 3 there exists a constant
$\;C>0\;$ such that 
$\;\;\sum\limits_{m\leq t}\varphi(m)\;\geq\;C\cdot t^2\;\;$
for all $\;t\geq 1\,$. (Actually, this estimate is certainly true
if we choose $\;C\,=\,{1\over 5}\,$.)
Hence we have
$$\sum_{(x,y)\in \G(T)}\varphi(x)\,\varphi(y)
\;\;=\;\;\sum_{1\leq y\leq T}\varphi(y)\cdot\sum_{1\leq x\leq T/y}\varphi(x)
\;\;\geq\;\;\sum_{1\leq y\leq T}\varphi(y)\cdot C{T^2\over y^2}\;.$$
Partial summation yields
$$\sum_{1\leq y\leq T}\varphi(y){1\over y^2}\;\;=\;\; 
{1\over T^2}\sum_{m\leq T}\varphi(m)
\;+\;\int\limits_1^T
{2\over u^3}
\Big(\sum_{m\leq u}\varphi(m)\Big)\;du\;\;\geq\;\;
C+2C\log T $$
and we arrive at (4.1) as requested.

\bp\mp
{\bff 5. Proof of Theorem 3}
\mp
Since the case $\;n=2\;$
is already settled by Theorem 2, in order to prove Theorem 3
we may assume $\;n\geq 3\,$. 
Further, by adapting the proof of Theorem 1 it is straightforward
to obtain the upper bound in Theorem 3.
Actually, this bound has the same order of magnitude
as the integral
$$I_n(T)\;=\;
\int\!\!\!\cdots\!\!\!\int\limits_{\!\!\!\!\!\!\!\!\!\!\!\!\D_n(T)} 
x_1\cdots x_n\;d(x_1,...,x_n)$$
with
$\;\,\D_n(T):=
\{\,(x_1,...,x_n)\in [1,\infty[^2\;|\;x_1\cdots x_n\leq T\,\}\;\,$
and it is plain to verify
\sp
\cl{$I_n(T)\;\asymp\;T^2(\log T)^{n-1}
\;\;\;(T\to\infty)$}
\sp
 for every $\;n\geq 2\;$                            
by induction starting from 
$\;I_2(T)\,=\,I(T,1,1)\;$ and using the
estimate $\;I(T,1,1)\,\asymp\,T^2\log T\;\,(T\to\infty)\;$ of Lemma 2.
\mp
On the other hand, following the lines of the proof of 
the lower estimate in Theorem 2 it is plain that
$$n!\cdot R_n^s(t)\;\;\geq\;\;
6^n\;\cdot\!\!\!\!\! \sum_{(x_1,...,x_n)\in\G_n(T)}
\!\!\!\!\!\varphi(x_1)\cdots\varphi(x_n)$$
with $\;T\,=\,n^{1-n}t\;$  and $\;\G_n(T)\,:=\,\D_n(T)\cap\Z^n\,$.
\mp
Now by applying Lemma 3 and partial summation, 
induction leads to
$$ T^2(\log T)^{n-1}\;\;\ll\;
\sum_{(x_1,...,x_n)\in\G_n(T)}
\!\!\!\!\!\varphi(x_1)\cdots\varphi(x_n)\;\;\;\;\,(T\to\infty)$$
for every $\;n\geq 2\;$ since
$$\int\limits_1^T
\(\Big(\log {T\over u}\Big)^{n-2}{2\log(T/u)\,+\,n-1\over u^3}\)\cdot u^2\;du
\;\;=\;\;{2\over n}\cdot(\log T)^n\,+\,(\log T)^{n-1}$$ 
for all $\;T\geq 1\;$ and every $\;n\geq 2\,$.
This concludes the proof of Theorem 3
\bp\mp
{\bff 6. Proof of Theorem 4}
\mp
The following facts, where always $\;k,h\in\Z\;$ is assumed,
are essential for our proof of Theorem 4.
\mp
(F1) {\it For every $\;k\geq 2\;$ and $\;h\geq 1\;$ there is at least one
irreducible $\;p(X)\in\Z[X]\;$ with $\;\deg p=k\;$ and $\;H(p)=h\,$.}
\sp
{\it Proof.} This is certainly true because, e.g., 
$\;X^k-hX^{k-1}-X^{k-2}-\cdots -1\;$ is irreducible,
which follows immediately from [4] Theorem 2.2.6.
\mp
(F2) {\it For every $\;k\geq 2\;$ and $\;h\geq 9\;$ the number of
all irreducible $\;p(X)\in\Z[X]\;$ with $\;\deg p=k\;$ and $\;H(p)=h\;$
is greater than $\;h^k/3\,$.}
\sp
{\it Proof.} We apply Eisenstein's Irreduciblity Criterion
with 2 as the testing prime. If $\,h\,$ is odd, then obviously
all polynomials $\;hX^k+2a_{k-1}X^{k-1}+\cdots+2a_1X+2\cdot(2l-1)\;$
with $\;l,a_i\in\Z\;$ and $\;2|a_i|<h\;$ and $\;4|l|<h\!-\!2\;$
are irreducible. If $\,h\,$ is even, then
all polynomials $\;(2l-1)X^k+hX^{k-1}+
2a_{k-2}X^{k-2}+\cdots+2a_1X+2\cdot(2l'-1)\;$
with $\;l,l',a_i\in\Z\;$ and $\;2|a_i|\leq h\;$ and $\;2|l|<h\;$
and $\;4|l'|\leq h\!-\!2\;$ are irreducible. 
Hence the requested number is not less than
$\;h^{k-1}(h\!-\!3)/2\;$ when $\,h\,$ is odd
and not less than
$\;(h+1)^{k-2}(h\!-\!1)(h\!-\!2)/2\;$ when $\,h\,$ is even.
\sp
Combining (F1) and (F2) we derive
\sp
(F3) {\it For every $\;k\geq 2\;$ and $\;h\geq 1\;$ the number of
all irreducible $\;p(X)\in\Z[X]\;$ with $\;\deg p=k\;$ and $\;H(p)=h\;$
is not smaller than $\;9^{-k}\cdot h^k\,$.}
\mp
On the other hand, since the number of
all $\;p(X)\in\Z[X]\;$ with $\;\deg p=k\;$ and $\;H(p)=h\;$
is certainly not greater than $\;2(k+1)(2h+1)^k\,$,
we have
\sp                      
(F4) {\it For every $\;k\geq 2\;$ and $\;h\geq 1\;$ the number of
all irreducible $\;p(X)\in\Z[X]\;$ with $\;\deg p=k\;$ and $\;H(p)=h\;$
is not greater than $\;2(k+1)3^k\cdot h^k\,$.} 
\mp
As usual, let us call a polynomial over $\Z$ {\it primitive} when
the greatest common divisor all of its coefficients is 1.
\mp
(F5) {\it For all $\;m\geq 1\;$ and $\;h\geq 1\;$ the total number of
all primitive polynomials $\,p(X)\,$ over $\Z$ with 
$\;\deg p=m\;$ and $\;H(p)=h\;$
is not greater than $\;2(m+1)3^m\cdot h^m\;$
and not smaller than $\;2^{m+1}\cdot\varphi(h)\cdot h^{m-1}\,.$} 
\sp
The upper bound corresponds to the bound in (F4) and is trivial. 
The lower bound comes from 
counting only all polynomials $\;\pm hX^m+aX^{m-1}+a_{m-2}X^{m-2}+\cdots+ a_0\;$
with $\;a,a_i\in\Z\;$ and $\;|a|,|a_i|\leq h\;$
where $\,h\,$ and $\,a\,$ are coprime.
\mp
Further, the following statement is obviously true. 
\mp
(F6) {\it If $\;1<{n\over 2}<k<n\;$ then for 
every $\;p(X)\,\in\,\Y_{k,n}(t)\;$ there
exists one and only one pair $\;(f,g)\in\Z[X]^2\;$
such that $\,g(X)\,$ is irreducible with 
$\;\deg g=k\;$ and $\,f(X)\,$ is primitive
and $\;f(X)\cdot g(X)\,=\,p(X)\,$.}
\bp
Now we a ready to prove Theorem 4.
Assume $\;1<{n\over 2}<k<n\,$. 
Then by (F6) the mapping 
$\;(f,g)\,\mapsto\,f(X)\cdot g(X)\;$
is a bijection from the set
\sp
\qquad$\big\{\,(f,g)\in\Z[X]^2\;\;\big|\;\;
\deg f\,=\,n-k\;\;\land\;\;\deg g\,=\,k\;\;\land\;\;H(fg)\;\leq\;t$

\rl{$\;\;\land\;\;(\,f\;\hbox{\it primitive}\,)\;\;\land\;\;
(\,g\;\hbox{\it irreducible}\,)\,\big\}\;$\qquad}
\sp
onto the set $\;\Y_{k,n}(t)\,$.
\bp
Consequently, with $\;t\ll T\ll t\,$, in view of Lemma 1 
and (F4) and (F5) we have
$$|\Y_{k,n}(t)|\;\;\ll\sum_{(x,y)\in \G(T)}x^{n-k}\cdot y^{k}
\;\;\ll\;\;T^{k+1}\;\;\;\;\;(T\to\infty)\leqno(6.1)$$
since $\;I(T;n-k,k)\ll T^{k+1}\;$ for $\;{n\over 2}<k<n\;$
by Lemma 2.
\bp
On the other hand, again with $\;t\ll T\ll t\,$, by Lemma 1 and 
by (F3) and (F5),
$$|\Y_{k,n}(t)|\;\;\gg\sum_{(x,y)\in \G(T)}\varphi(x)x^{n-k-1}\cdot y^{k}
\;\;\;\;\;(T\to\infty)\;.$$
By writing
$$\sum_{(x,y)\in\G(T)}\varphi(x)x^{n-k-1}\cdot y^{k}\;\;=\;\;
\sum_{1\leq x\leq T}\varphi(x)x^{n-k-1}
\sum_{1\leq y\leq T/x} y^{k}$$
and applying the trivial estimate
\sp
\cl{$\sum\limits_{1\leq y\leq u} y^{k}\;\;\geq\;\;\int\limits_0^u y^kdy\;-\;u^k
\;\;=\;\;{1\over k+1}u^{k+1}\,-\,u^k\;\;(u\geq 1)$}
\sp
and Lemma 4 with $\;\alpha\,=\,n-2k-2\,<\,-2\;$ on the one hand and
with $\;\alpha\,=\,n-2k-1\,<\,-1\;$ on the other, we derive
the desired lower estimate
$$T^{k+1}\;\;\ll\sum_{(x,y)\in \G(T)}\varphi(x)x^{n-k-1}\cdot y^{k}
\;\;\;\;\;(T\to\infty)\;.$$
\mp
Further, the estimate
$\;\;t^2(\log t)^2\,\ll\,|\,\Y_{3}(t)\setminus\Y_{2,3}(t)\,|\,\ll\,
t^2(\log t)^2\;\;$
is equivalent to Theorem 3 for $\;n=3\,$.
In particular, the factor $\,(\log t)^2\,$ in the estimate 
in Theorem 4 cannot be omitted when $\;n=3\,$.
\mp
In order to verify $\;\;|\,\Y_{n}(t)\setminus\Y_{n-1,n}(t)\,|\,\ll\,
t^{n-1}\;\;$ for $\;n\geq 4\;$ we note that                         
$$\Y_{n}(t)\setminus\Y_{n-1,n}(t)\;\;\,\subset\;\;\,
\Y_n^*(t)\,\cup\!\!\!\!\!\!\!\!\!
\bigcup\limits_{\;\;\;{n\over 2}<k\leq n-2}
\!\!\!\!\!\!\!\!\!\!\Y_{k,n}(t)\leqno(6.2)$$
where 
$\,\Y_n^*(t)\,$ is the set of all reducible polynomials $\,p(X)\,$
over $\Z$ with degree $\,n\,$ and height $\;\leq t\;$
such that the degree of every irreducible factor
of $\,p(X)\,$ is not greater than $\,{n\over 2}\,$.
\sp
Now, for every $\;p\in\Y_n^*(t)\;$ 
we can write $\;p(X)\,=\,f(X)\cdot g(X)\;$ 
with $\;f(X),g(X)\,\in\,\Z[X]\;$ such that
the degrees of $\,f(X)\,$ and $\,g(X)\,$
are both not greater than $\;n-2\,$.
Hence, by following the arguments in Section 3 we only have
to estimate the sum in (3.2) for $\;2\leq k\leq n\!-\!2\;$
in order to obtain $\;|\Y_n^*(t)|\,\ll\,t^{n-1}\,$.
Thus, in view of (6.1) via (6.2) we arrive at
$\;\;|\,\Y_{n}(t)\setminus\Y_{n-1,n}(t)\,|\,\ll\,
t^{n-1}\;\;$ for $\;n\geq 4\;$ and
this concludes the proof of \hbox{Theorem 4.}
\bp\bp
{\bf Final Remark.} In view of our proofs it is not difficult
to find explicit bounds $\;C_n\;$
in Theorem 1 and to  
produce explicit $\ll$-constants 
for all estimations in Theorems 2, 3, 4 
which depend only (and in a simple way) on the degree 
$n$.

\bp\bp\bp
{\bff References}
\bp
[1] K.~{\smc Chandrasekharan}:
{\it Introduction to Analytic Number Theory}. 
Springer 1968.
\mp 
[2] K.~{\smc D\"orge}: {\it Absch\"atzung der Anzahl der reduziblen Polynome.}
Math.~Ann.~{\bf 160}, 59-63 

\rl{(1965).}

[3] G.~{\smc P\'olya} and G.~{\smc Szeg\"o}:
{\it Problems and Theorems in Analysis}$\,$ Vol.~II. Springer 1976.
\mp
[4] V.V.~{\smc Prasolov}: {\it Polynomials.}
Springer 2004.

\bp\mp\bp
{\bff Author's address}
\mp
Institute of Mathematics 
\sp 
University of Natural Resources and Life Sciences
\sp
Vienna, Austria
\sp
{\sl E-mail:} {\tt gerald.kuba@boku.ac.at}
\bp\bp\mp\hrule\bp\bp\mp
{\sl This paper has been published in {\bf Mathematica Slovaca 59 (2009)}.}

\end